\documentclass[12pt]{article}
\usepackage{amsmath,amsfonts,amsthm,amssymb,amscd}
\usepackage[all]{xy}
 \usepackage[oldstylenums]{kpfonts}
   \usepackage[cal=boondoxo]{mathalfa} 
    \usepackage{hyperref} 

\newtheorem{thm}{Theorem}[section]

\newtheorem*{thm*}{Theorem}
\newtheorem*{corr*}{Corollary}
\newtheorem{lemma}[thm]{Lemma}
\newtheorem{prop}[thm]{Proposition}
\newtheorem*{prop*}{Proposition}

\newtheorem{corr}[thm]{Corollary}
\theoremstyle{definition}
\newtheorem{dfn}[thm]{Definition}
\newtheorem{exmple}[thm]{Example}
\newtheorem{exmples}[thm]{Examples}

\theoremstyle{remark}
 
\newtheorem*{rmq}{\textit{Remark}}
\newtheorem{rmk}[thm]{\textit{Remark}}

\def\C{{\mathbf{C}}}
\def\G{\mathbf G}
\newcommand{\Q}{{\mathbf{Q}}}
\def\bP{\mathbf{P}}
\newcommand{\R}{{\mathbf{R}}}

\def\bU{\mathbf U}
\def\bV{{\mathbf V}}
\def\bW{\mathbf W}
\newcommand{\Z}{{\mathbf{Z}}}

\def\eh{\mathfrak{h}}
\def\ek{\mathfrak{k}}
\def\ep{\mathfrak{p}}
\def\eg{\mathfrak{g}}
\def\AA{\mathcal{A}}

\newcommand\CC{{\mathcal C}} 
 
\newcommand\FF{{\mathcal F}} 
\newcommand\GG{{\mathcal G}}

\newcommand\OO{{\mathcal O}} 
\newcommand\XX{{\mathcal X}} 
\newcommand\VV{{\mathcal V}} 
 
\def\ii{{\rm i}}

\newcommand{\comp}{\raise1pt\hbox{{$\scriptscriptstyle\circ$}}}
\def\del{\partial}
\def\id{\mathop{\rm id}\nolimits}
\def\rank{\mathop{\rm rank}}
\renewcommand\setminus{-}
\def\into{\hookrightarrow}

 \def\lset{\{}  
\def\rset{\}}  
\def\set#1{\lset#1\rset} 
  
\def\mapright#1{\mathop{\vbox{\ialign{
                                ##\crcr
    ${\scriptstyle\hfil\;\;#1\;\;\hfil}$\crcr
 \noalign{\kern2pt\nointerlineskip}
    \rightarrowfill\crcr}}\;}}

\def\mapleft#1{\mathop{\vbox{\ialign{
                                ##\crcr
    ${\scriptstyle\hfil\;\;#1\;\;\hfil}$\crcr
 \noalign{\kern2pt\nointerlineskip}
    \leftarrowfill\crcr}}\;}}

\newcommand\gl[1]{\operatorname{GL}({#1})}
\newcommand\so[2]{\operatorname{SO}^{#2}({#1})}
\newcommand\su[1]{\operatorname{SU}({#1})}
\newcommand\ogr[2]{\operatorname{O}^{#2}({#1})}
\newcommand\ugr[2]{\operatorname{U}^{#2}({#1})}
\newcommand\slgr[1]{\operatorname{SL}({#1})}
\newcommand\smpl[1]{\operatorname{Sp}({#1})}

\newcommand\gr{\operatorname{Gr}}

\newcommand\Hom{\mathop{\rm Hom}\nolimits}
\newcommand\End{\mathop{\rm End}\nolimits}
\newcommand\aut{\operatorname{Aut}}
\def\half{\frac{1}{2}}

\def\aut{\mathop{\rm Aut}}
\def\Res{\mathop{\rm Res}}
\def\bS{\mathbf S}
\def\Lie{\mathop{\rm Lie}}
\def\mor{\mathop{\rm Mor}}
\title{Rigidity of Spreadings and Fields of Definition}
\author{Chris Peters\\
Technical University Eindhoven \\
and Universit\'e Grenoble Alpes \\
email: {\tt c.a.m.peters@tue.nl}  
}
\date{November 2016}

\begin{document}
\maketitle
\begin{abstract} Varieties without deformations are  defined over a  number field. Several old and new examples of this phenomenon are discussed such
as  Bely\u \i\ curves and Shimura varieties. Rigidity is related to maximal Higgs fields which come from variations of Hodge structure.  Basic properties 
for these due  to P. Griffiths, W. Schmid, C. Simpson and, on the arithmetic side, to Y. Andr\'e and I. Satake all play a role.
This note tries to give a largely self-contained exposition of these manifold ideas and techniques, presenting, where possible, short new proofs for  
key results.

\centerline{AMS Classification {\tt 11G18, 14C30, 14D07, 14G35} }
\end{abstract}

\section*{Introduction}

Results stating that certain types of algebraic  varieties are definable over a number field are
scattered in  the literature. Arguably, those    most studied form  the class of Shimura varieties \cite{nodefs,arrigid,milne}.
Another  famous example is Bely\u \i's theorem \cite{3points} which characterizes curves over $\overline\Q$
as those which have a  Bely\u \i  \ representation, i.e.,  a branched cover of the line branched in exactly
three points.  In dimension two,  the fake projective planes  \cite{fake1,fake2}
and  the Beauville surfaces \cite[Exercise X.13 (4)]{beau}, \cite{beauville2} are known to have models over $\overline\Q$.

Such examples can be uniformly explained by   constructing a   suitable spread  of the varieties concerned as demonstrated in Sections~\ref{sec:Basic} and \ref{sec:exmples}.
 
Of a totally different flavor are the applications to special subvarieties of  Shimura varieties in  Section~\ref{ssec:SS}
on the one hand, and to  splittings of Higgs bundles as given in Section~\ref{ssec:rigidss} on the other hand. As has been known since the work of
Viehweg and Zuo \cite{carshcurves}, the last two are just facets of the same phenomenon:  Higgs bundles of a very special kind, those that they  called
"maximal" are directly related to special subvarieties of certain Shimura varieties. One of the goals of this note is to
show that rigidity plays a central role in this; exploiting this, simplifies several of the arguments.

The overall goal of this note  is to  show how a few relatively simple ideas  plus some standard  techniques from deformation theory and Hodge
theory explain  a   wide range of  phenomena of the above kind.
It brings together various known results from very different subfields of mathematics. This is the reason why   I thought  to
 explain some of the basic notions and techniques from these fields,  and also   to search for new simpler proofs.

This note  has been inspired by discussions with Stefan M\"uller-Stach.  Equally influentual has been    \cite[Ch. 4]{tangac} as well as the last chapter of \cite{redbook}. 
Finally I would  like to thank Christopher Deninger for pointing out the references \cite{arrigid,nodefs}.

\newpage 
\tableofcontents

\section{Spreads of  varieties and rigidity}
\label{sec:spreadrigid}

\subsection{Spreads}
\label{ssec:spread}

The  "spread philosophy"  roughly states that a complex algebraic variety can be seen as a family 
 over a  base variety  determined by specifying some
 transcendence  basis of  the   field of definition of the  variety. Spreads  are by no means unique
 but all share the crucial property that,  by construction,  the total family is always defined over a number field.
 
 Although this construction can be phrased in
varying  generality \cite[\S 4.1]{tangac}, the following somewhat restricted version suffices for this note.
\begin{prop}\label{spread} 
Let $X$ be  smooth complex quasi projective variety.  
There exists a  smooth   family  
$f:\XX \to S$ defined over $\overline\Q$ such that
\begin{enumerate}
\item $\XX, S$ are smooth quasi-projective;
\item there is a canonical point $o\in S$ such that $f^{-1}o=X$;
\item if $s\in S$ is a $\overline\Q$--rational point, the fiber $X_s=f^{-1}s$ is  defined over $\overline\Q$.
\end{enumerate} 
\end{prop}
\proof
Suppose for simplicity that the variety  $X$ is projective and is given by a finite set of polynomials.  The coefficients of these polynomials
generate  a  field $k$ of finite transcendence degree $r$
over $\Q$, say $k=k_0(\alpha_1,\dots,\alpha_r)$ where $\set{\alpha_1,\dots,\alpha_r}$ is a transcendence basis  for $k$  and where $k_0$ is a number field, say of the form $\Q[x_0]/P$ with $P$ some monic irreducible polynomial. Then $k$ is  the function field of some    complex algebraic variety $S'$. The deformation will be constructed over a Zariski open subset $S$ of $S'$.

The basic idea is to replace the coefficients $\alpha_j$  of each of the  polynomials defining the variety $X$, by  variables $x_j$.
A  point $s\in S'$ corresponds to a field  $k(s)$ isomorphic to $k$.  If one replaces the coefficients in $k$ of  a defining set
of  homogeneous  equations   for $X$ by the corresponding coefficients in $k(s)$ one gets a variety $X_s$. 
The $X_s$ glue to a variety   $\mathcal{X}$ fibered over $V(P)$. Indeed, it is given by the same equations as $X$ except that
the coefficients for these equations are not considered as numbers but as   $\Q$--polynomials in the supplementary variables $x_j$ tied by the extra equation $P(x_0,\dots,x_r)=0$.

Substituting $x_j=\alpha_j$ gives a canonical $k$--valued point $o\in S'$ and
by construction $X_o=X$. 
Since $k/k_0$ is separable, this point is a non-singular point. Now replace $S'$ by a suitable Zariski open neighborhood $S$ of $o$ such that the variety $S$ is smooth.  Again by separability,   this variety  is smooth along $X_o$.
But it  might still be singular or reducible. To remedy this, first take the 
component of $\mathcal{X}$ which contains $f^{-1}o$. Then replace $S$ by a smaller neighborhood of $o$ such
that not only the total space is smooth, but also all of the fibers of the fibration are smooth.
The resulting  family, still denoted  $f:\XX \to S$,  is a smooth deformation of $X$.
By  construction, the Zariski-open subsets  figuring in the  construction  are complements of equations  over $\overline \Q$, and so the resulting  family is defined over $\overline \Q$.

Finally, since $\overline \Q$ is algebraically closed, $S$ contains points $s$ defined over $\overline \Q$.
This amounts to replacing the variables $x_j$  figuring in the coefficients for the equations of $\mathcal X$   by suitable algebraic numbers and hence $X_s$ is defined over $\overline\Q$.
\endproof
\begin{rmk} \label{morespreads}
There are several variants of this result: one can spread pairs $(X,Z)$ with $Z$ a closed subvariety of $X$.
Similarly, one can spread a given morphism $f:X\to Y$ between varieties.
\end{rmk}

\subsection{Deformations and rigidity}
\label{ssec:rigid}

Let me first recall some basic definitions and facts. More details and proofs can be found for example in  \cite{sernesi}.
\subsubsection*{Kodaira-Spencer classes}

A complex  variety $X$ is said to be  infinitesimally, respectively locally  \emph{rigid}  if any infinitesimal deformation
of $X$, resp. any local deformation $p:\XX\to S$ of $X$ with $S$ sufficiently small,  is trivial, i.e. isomorphic to the product deformation.
This can be rephrased by saying that  if, say   $o\in S$ is such that the fiber of $p$ over it is isomorphic to $X$, say $\iota :X_o=f^{-1}o\mapright{\simeq} X$,  then   there is a morphism
 $S\to \aut(X),\quad s\mapsto g_s,\, g_o=\id_X$
 inducing a product structure on  the family $\XX \to S$:
 \[
X\times S \mapright{\simeq}  \XX, \quad (x,s) \mapsto (\iota \comp g_s(x),s).
\]
As is well known, a  variety $X$ is indeed  infinitesimally or locally rigid if  $H^1(X,\Theta_X)=0$. If such a variety appears in a deformation
$p:\XX\to S$ of $X\simeq X_o$, $o\in S$, finer information is present by looking at  the
 the \emph{Kodaira-Spencer map}. Recall that it is given as  the extension class of the exact sequence
\begin{equation}
\label{eqn:KSasExtclass}
0 \to T_o(S)\otimes \OO_X \to \Theta_{\XX}|X \to \Theta_X\to 0
\end{equation}
of $\OO_X$--modules. In other words, it gives a characteristic map
\[
\kappa_p : T_o(S)  \to H^1(X,\Theta_X).
\]
For a given deformation, it   measures deviation of triviality of the deformation:
\begin{thm}[\protect{\cite[Thm. 18.3]{KS}}] \label{ks}
Suppose that a  family   $p:\XX\to S$   is regular in the sense that $\dim H^1(X_s,\Theta_{X_s})$ is constant for  $s\in S$. Then  it  is trivial if and
only if  $\kappa_p=0$. 
\end{thm}
Observe that this theorem gives back the criterion that $X$ is rigid if and only if   $H^1(X,\Theta_X)=0$.
Indeed, if this is the case, by semi-continuity, any sufficently small deformation of $X$ is regular and the theorem applies to show
rigidity.


\subsubsection*{Variants}
1. Infinitesimal  \emph{deformations of 
 pairs} $(X,Z)$ with $Z$ a closed subscheme of a \emph{smooth} variety $X$.  Any such deformation $p$ with base $(S,o)$ (i.e. with 
  fiber over $o$  isomorphic to $(X,Z)$) is classified by its  Kodaira-Spencer map
  \begin{equation}
  \label{eqn:KSpairs}
  \kappa_p:  T_oS \to H^1(X, \Theta_X(Z)),
  \end{equation}  
   where $\Theta_X(Z)$ is the sheaf of germs of vector fields on $X$ tangent to $Z$\footnote{This is Sernesi's notation; if $Z$ is a normal crossing divisor
   it is dual to $\Omega^1_X(\log Z)$ and other auhors use $\Theta_X(-\log Z)$ in this case.}. This deformation is rigid precisely
 when $\kappa_p=0$ as before.
\\ 
2.  \emph{Deformations of morphisms}  $f:X\to Y$. These are given by a commutative  diagram
 \[
 \xymatrix{ X  \ar[rr]_{ f}    &&  Y \\
  X _o   \ar[u]_{\iota}^{\simeq}  \ar@{^(->}[d]   \ar[rr]_{ F|X_o} &&  Y_o \ar[u]^{\iota'}_{\simeq} \ar@{^(->}[d]\\
  \mathcal  X\ar[rr]_F \ar[dr]_{p_1} && \mathcal Y \ar[dl]^{p_2} \\
 & S.
 }
 \]
 A deformation of morphism as above  is a \emph{deformation keeping the source, respectively  target fixed} if $p_1$ resp. 
 $p_2$ are a trivial deformations.  
  A  morphism  $f$  is rigid, if all infinitesimal deformations of $f$
  are trivial   in the sense that there are  morphisms
 \[
 \aligned
S\to \aut(X)&\quad   s\mapsto g_s,\, g_o=\id_X\\
S\to \aut(Y)& \quad  s\mapsto g'_s,\, g'_o=\id_Y.
\endaligned
 \]
 which trivialize  the deformation: for all $s\in S$ there is a commutative diagram
\[
 \xymatrix{
 X  \ar[d]_{\simeq}^{\iota\comp g_s}   \ar[rr]^{f} && Y \ar[d]^{\simeq}_{\iota'\comp g'_s} \\
X_s \ar[rr]_{F|X_s} && Y_s.
 }
 \]

  Two special cases will be used in this note:
 \begin{enumerate}
 \renewcommand{\labelenumi}{\alph{enumi}) }
  \item Deformations of a morphism $f:X\to Y$ between non-singular varieties keeping source and target fixed.
Such morphisms are classified by  the vector space $H^0(X,f^*\Theta_Y)$.
\item  Deformations of closed embeddings $f: Z\into X$ between smooth varieties with target fixed.
 Here the characteristic morphism is
 \[
 \kappa_F :  T_oS \to H^0(Z, N_{Z|X}),
 \]
where $N_{Z|X}$ is the normal bundle of $Z$ in the ambient manifold $X$. Note that  automorphisms of $X$ 
yield non-trivial deformations of $f$ but these  are trivial as deformations of $Z$ itself. Indeed, there is an exact sequence
\[
0\to H^0(Z,\Theta_Z)  \mapright{i^*}  H^0(Z,\Theta_X|Z ) \to  H^0(Z, N_{Z|X}) \mapright{\delta} H^1(Z,\Theta_Z).
\]
The quotient $H^0(\Theta_X|Z )/i^*H^0(\Theta_Z)$ is the space of isomorphisms classes of infinitesimal deformations of $f$ keeping $Z$ and $X$ fixed; the next term in
the sequence, $ H^0(X, N_{Z|X})$,   is the space of infinitesimal deformations of $f$ keeping only $X$ fixed and $\delta$ maps such a deformation to the corresponding deformation
of $Z$, i.e., it is  the forgetful map.
The embedding $f$ is rigid in this case precisely if $\kappa_F=0$.  If $Z$  itself is rigid,  this would  follow if   $H^0(Z, \Theta_Z )\to H^0(Z, \Theta_X|_Z )$ is  surjective.
In case $Z$ admits no vector fields, this means  $H^0(Z, \Theta_X|_Z )=0$.  
For later use, consider the following special case,  that of a totally geodesic submanifold $Z$ of $X$:
\begin{prop}[\protect{\cite[Sect. 11.5]{3authors}  }]\label{GeodEmbedded} Let $X$ be a manifold  equipped with a hermitian metric, and let $Z\subset X$ be a totally geodesic submanifold for which $H^0(Z,\Theta_Z)=0$. Then 
the tangent bundle sequence for $Z\subset X$ splits. Hence $ H^0(Z,\Theta_X|_Z)=0\iff H^0(Z, N_{Z|X})=0 $.
It follows that  $Z$ is rigidly embedded (keeping the target fixed) if and only if  the embedding is  rigid keeping source and target fixed.
\end{prop}
In particular,  since $\delta$ is the zero map in this case, it is irrelevant whether  $Z$ itself is rigid or not.
\end{enumerate}

\subsubsection*{Kodaira-Spencer classes and  spreading}  
The Kodaira-Spencer class of the spread  family $f:\XX \to S$   from  Prop.~\ref{spread} incorporates arithmetic information, since the dual of $T_o(S)$ is
the  complex vector space $\Omega_{k/\Q}\otimes_\Q\C$. Also, $\Omega^1_{\XX}|X= \Omega^1_{X/k}$, the sheaf of  K\"ahler differentials  on the $k$-variety $X$.
The dual of the exact sequence \eqref{eqn:KSasExtclass} then reads 
\[
0 \to  \Omega^1_{X} \to \Omega^1_{X/k} \to  \Omega_{k/\Q}\otimes_\C \OO_X     \to 0.
\]
The extension class of the dual of the above  sequence is  the Kodaira-Spencer class
for the spread family $f:\ \XX \to S$ at $o$. It depends on the choice of the field $k$:
\begin{equation}\label{eqn:kspcl}
\kappa_{X/k}  \in H^1(\Hom_{\OO_X}(\Omega^1_{X}  ,\Omega_{k/\Q}\otimes_\C \OO_X))\simeq  \Hom_\C(T_o S, H^1(X,\Theta_X)).
\end{equation}

\begin{corr} \label{spreadregular} The spread family from Prop.~\ref{spread} is regular. It is a trivial  deformation  if  and only  if the Kodaira-Spencer class \eqref{eqn:kspcl} vanishes.
\end{corr}
\proof  
 To see  regularity, first observe that $\dim H^1(X,\Theta_X)$ depends on the isomorphism class of $X$ as an abstract  algebraic variety. Secondly, since 
  all $X_s$, $s\in S$ with the property that $s$
 corresponds to a transcendental number are mutually isomorphic as abstract     algebraic varieties,  $\dim H^1(X_s,\Theta_{X_s})$
 is the same for all such $s\in S$.  This set corresponds to points in $S$  not lying on any
 proper subvariety of  $S$ and hence is dense in $S$. Upper semicontinuity 
 of  $\dim H^1(X_s, \Theta_{X_{s}})$  then implies that this dimension  is locally constant, i.e., the family is regular.   The
 result follows from  Theorem~\ref{ks}.
\endproof

\subsection{Rigidity and fields of definition}
\label{sec:Basic}
 
\begin{prop}\label{model}
{\rm 1)} Let  $X$ be a smooth complex quasi-projective variety.  
Assume that the Kodaira-Spencer class \eqref{eqn:kspcl}  of some spread family of $X$ vanishes (e.g. in case  $X$ is rigid). Then $X$ has a \emph{model  over a number field}, i.e.,
\[
X  \simeq X' \otimes _{\overline\Q} \C, \quad X'  \text{ is defined over $\overline\Q$,}
\]
and where the isomorphism is  defined over $\C$. This model is unique if $H^0(\Theta_X)=0$.
\\
{\rm 2)} Let $(X,Z)$ be a pair of varieties, where $X$ is smooth and $Z\subset X$ a closed embedding. Assume that the Kodaira-Spencer class \eqref{eqn:KSpairs} of a  spread family
for $(X,Z)$ vanishes (e.g. in case $(X,Z)$ is rigid).  Then the pair  $(X,Z)$ has a model over a number field, The model is unique if  $H^0(\Theta_X(Z))=0$.
\\
{\rm 3)}  In the relative situation of a morphism $f:X\to Y$ between complex quasi-projective varieties suppose that $f$ is rigid. Then $X,Y$ and $f$ have a model over $\overline \Q$.\\
{\rm 4)} In the relative situation, suppose that $Y$ is defined over $\overline\Q$ and that $f$ is rigid fixing the target. Then the same conclusion as in {\rm 3)} holds.
\end{prop}
\proof
1)  Rigidity implies that the fibers of any sufficiently small deformation of $X$  are isomorphic to $X$.
This holds in particular  for the spread   $f:\mathcal{X}\to S$  from Prop.~\ref{spread}.
So, if  $s\in S(\overline \Q)$,  one has an isomorphism  $X_s\simeq X_o=X$   and  since $X_s$ is defined over $\overline \Q$,
$X$ has a model over $\overline \Q$. If, moreover, $H^0(\Theta_X)=0$ there is no non-trivial deformation of $\id_X$ and the isomorphism $X_s\simeq X_o$ is unique (compare with the definition above).
\\
2)  The argument is as for 1), using an obvious variant of  Prop.~\ref{spread} for pairs. See remark~\ref{morespreads}.

Note that  3) and 4) 
can be reduced to embeddings, since $f$ is rigid if and only the embedding of graph of $f$ in $X\times Y$ is a rigid morphism, 
and the graph is defined  over $\overline\Q$ precisely when $f$ is.
For embeddings $i:X\into Y$, to find a variety over which to spread, start with equations for $Y$ and let $k_1$ be the field extension of $\Q$ obtained by adjoining the coefficients. The embedding is then specified by supplementary equations whose coefficients are adjoined to $k_1$. The resulting field $k=\overline \Q(S)$  is the function field of the base variety  $S$. Observe that if the variety $Y$ is defined over a number field, $k_1$ is also a number field and then $S$ parametrizes a deformation of $X$ in the fixed variety $Y$.
Rigidity in both cases ensures that the embedding has a model over a number field.
 \endproof

\begin{exmples}

1. \emph{Fake projective planes} are compact complex surfaces of general type with $p_g=q=0$ and with 
$K^2=9$. They are known to be quotients of the complex unit $2$--ball by  an arithmetic subgroup, and are
also known to be rigid. See \cite{fake1,fake2}.
\\
2. Let $S$ be a \emph{Beauville surface} \cite[Exercise X.13.(4)]{beau} and \cite{beauville2}. These are certain  minimal surfaces of general type with $K^2=8$, $p_g=q=0$.
Such a surface is rigid \cite{fibres} and so, by Proposition~\ref{model}, it has a model over $\overline\Q$. Its complex
conjugate cousin, also a Beauville surface,  is   rigid as well.

By \cite{beauville2}, there are a  two more types of surfaces similar to Beauville's examples in that they are all quotients of
a product of two curves of genera $>1$ by a freely acting finite group $G$ and having moduli spaces of dimension $0$.  Here $G$ is one of two non-abelian  groups of order 256. The first 
gives an example whose moduli space consists of   three
$0$--dimensional components, the second group leads to  a unique example.

\end{exmples}

The next result gives an application   in the relative setting. It leads up to Bely\v \i\ curves:

 \begin{prop} \label{belyi} Suppose $X,Y$ are smooth projective of the same dimension, $p:X\to Y$ is a surjective
finite morphism  with smooth branch locus $B\subset Y$. Assume that
  $Y$ is rigid and that $B$ is rigidly embedded in $Y$.
Then $X$ has a model  over a number field. 
\end{prop}
\proof
One constructs a  spread of the  morphism $p:X\to Y$ as in the proof of Cor.~\ref{model}.3. Call it $\tilde p: \XX\to Y\times S$.
We do not now that $p$ is rigid. But the induced
deformation of $f$, the family $\XX \to Y\times S\to S$, is differentiably  locally trivial over $S$  and so the topological structure
of the fibers $p_s$ of the map $\tilde p$  does not vary.   Away from the branch locus,
the map $p_s $ is a finite \'etale cover and so the complex structure
on 
$$
X_s^0 :=p^{-1}_s(Y \setminus  B_s ) \subset X_s 
$$ 
 is locally determined by the complex structure on $Y \setminus  B_s $, which by rigidity of the embedding of $B$ in $Y$ is 
 independent of $s$. The manifold  structure of $X_s$ is fixed and so it only has to be checked that the complex structure on
 it is completely determined by the complex structure on the Zariski open subset $X_s^0 $.
 
 To show this, note that holomorphic functions on 
$Y $ are bounded near the branch locus and so, by 
Riemann's extension theorem, their lifts to $X_s^0 $  can be extended uniquely to  $X_s$. So indeed,   up to isomorphism,   the
complex structure on  $X_s $ does not depend on $s$.   
As before, pick any $s\in S $ defined over $\overline  \Q$
(which exists since $S$ is by construction defined over $\overline \Q$). Then, not only $Y_s$ is defined over $\overline  \Q$, but also
$X_s$ is, and hence, by rigidity,  so is  the variety  $X_o=X$.
\endproof

\begin{rmq}
By \cite[p. 468--473]{redbook}, a variant of the above  proof is apparently due to  Carlos Simpson.
\end{rmq}

\begin{exmples}
1. Recall that a  \emph{Bely\u\i \ curve} \cite{3points} is a complex projective curve admitting a cover to $\bP^1$
ramified only in the three points $0,1,\infty$. Three distincts  points in $\bP^1$ define a rigid divisor since
three distinct points can always be mapped to three given distinct points by a projective transformation of $\bP^1$. 
Bely\u\i\ showed (loc. cit.) that a complex projective curve can be defined over $\overline\Q$
if and only if it is isomorphic to a Bely\u\i\ curve. The above Proposition shows that the fact
that Bely\u\i \ curves are defined over a number field is an example of a quite general phenomenon.
The converse statement however requires an explicit construction which is very particular to curves.
See \cite[Sect. 9.2]{redbook} for a proof in the style of this paper.
%
\\
2. For higher dimensional examples, including branched covers of $\bP^2$ branched
in $4$ or less lines,  see \cite{kapil}.
\end{exmples}

\section{Further examples of models over number fields}
\label{sec:exmples}

%

 \subsection{Locally symmetric spaces}
\label{ssec:locsymspace}

Let $D=G(\R)/K$ be a hermitian symmetric domain,  $\Gamma$  a torsion free arithmetic subgroup of  $G(\R)$ and let $X=\Gamma\backslash D$
be the corresponding locally symmetric space.  Such $X$ give examples of Shimura varieties  for which it is known
 that they can be defined over a number field.   See e.g.  \cite{milne} for background.
 Shimura varieties will be investigated more in detail below in Section~\ref{ssec:shimvars}.

Here I want to present another approach, due to Faltings which is more in the spirit of this note.   
\begin{prop}[\protect{\cite{arrigid}}]  \label{faltings} The pair 
 $(\overline X,\del X)$ has a unique model over $\overline \Q$.
\end{prop}
\proof
 I give a sketch of Faltings' proof.\footnote{For more details and a generalization of the results of \cite{cv} to the non-compact situation see \cite{rigidlocsym}.}
 The specific Kodaira-Spencer class  $\kappa_{(\overline X,\del X)}$
 coming from the  derivations of $\C/\Q$  
given by  \eqref{eqn:kspcl}  lands  in  the vector space 
$H^1(\overline  X,\Theta_X(\del X))$ measuring infinitesimal  deformations of the pair $(\overline  X,\del X)$. 
Using harmonic theory, Faltings shows  that each of these can be  
  represented by  a unique vector valued harmonic form  $\mathsf{H} \kappa_{(\overline  X,\del X)}$  on $D$ of type $(0,1)$. Moreover,
  the assigment $(\overline X,\del X)\mapsto  \mathsf{H}  \kappa_{(\overline X,\del X)}$ is functorial and equivariant with respect to group actions.

 Using this property for the various Hecke correspondences, one shows that such a  harmonic  form is
  $\Gamma$--invariant for all possible arithmetic subgroups $\Gamma\subset G$. This form lifts to   $D$ as a
 $G(\Q)$--invariant  harmonic $1$-form  with values in the tangent bundle. 
By density  it is then $G(\R)$--invariant on $D$. But such a form must vanish.  One sees this as follows.
By  \cite[Ch. VIII, \S 7]{Helg2} the complex structure on the tangent space $T_oD$ at a any point $o$ of the  hermitian symmetric domain $D$ 
is induced from the action of the center $Z \simeq \mathbf{U}^1$ of the isotropy group of $o$  on $D$: 
 $z\in Z$ induces  multiplication with $z$. Hence, if   $\alpha$ is  a global $(0,1)$-form on $D$ with values in the tangent bundle
  point $o$  the action is given by $z^*(\alpha)= (\bar z^{-1}\cdot z) \cdot \alpha$.  So $Z$-invariance, implies  $\alpha(o)=0$. Since $o$ is arbitrary,
$\alpha=0$.
 
 Next, one observes that the spread family for the pair $(\overline  X,\del X)$
 is regular  in the Kodaira-Spencer sense. The proof is similar to the proof of Cor.~\ref{spreadregular}.
 Hence one may apply  (a relative variant of) Theorem~\ref{ks}:  $(\overline X,\del X)$ is rigid, and hence this pair has  model over $\overline \Q$.
 
 Uniqueness then follows from $H^0(\overline X,\Theta_{\overline X}( \del X))=0$ (no vectorfields can be tangent along the boundary divisor). Faltings gives an explicit argument reducing the proof to the assertion that there exists  no $G(\R)$--invariant holomorphic
 vector fields on $D$. For the last assertion in loc. cit. no proof is given, but the argument is similar  to what we did before: The element $z\in Z=$ \{center of the isotropy
 group of $G(\R)$ at $o$\} acts as multiplication with $z$ on tangent vectors at $o$ and so, invariance implies that
 any  global tangent vector field on $D$  invariant under the action of $G(\R)$ vanishes at $o$ and hence everywhere.
 \endproof

\subsection{Holomorphic maps into     locally   symmetric spaces }
\label{ssec:intolocsymsp}

 As before, let   $X=\Gamma\backslash G/K$
be a locally symmetric space  of hermitian type. 
To $D=G/K$ and a parabolic subgroup $P\subset G$ one associates a boundary component $D(P)$
which is also a bounded symmetric domain. Introduce
\[
\text{rank of $D$}= \ell(D) = \min_P \dim D(P).
\]
The numbers $\ell(D)$ for $D$ irreducible are collected in Table~\ref{hsds}.
\begin{table}[ht]
\caption{Hermitian symmetric domains}
\begin{center}
\begin{tabular}{|c|c| c|}
\hline
Domain & $\dim D$ & $\ell (D)$\\
\hline
$I_{p,q} =\su {p,q}/\text{S}(\ugr {p}{} \times \ugr q {} )$& $pq $& $(p-1)(q-1)$\\
\hline
$II_g= \so {2g} {*}/\ugr g{}$ & $\half g(g-1)$ & $\half (g-2)(g-3)$\\
\hline
$III_g=\smpl{g}{}/\ugr g{}$& $\half g(g+1)$  & $\half g(g-1)$ \\
\hline
$IV_n=\so {2,n}{o}/ \so{2}{}\times \so {n}{}$ & $n$ & $1$  \\
\hline
$V= E_6/\so {10}{}\cdot \so 2 {}$& $16$ & $1$\\
\hline
$VI=E_7/ E_6 \cdot \so 2 {}$& $27$ & $8$\\
\hline
\end{tabular}
\end{center}
\label{hsds}
\end{table}%
The rigidity result I use here is due to Sunada:
\begin{prop}[\cite{sun}] With the above notation, let $M$ be  projective, $f: M\into X=\Gamma\backslash D$  with $X$ compact,  is rigid 
keeping source and target fixed, whenever $\dim M \ge \ell(D)+1$.
\end{prop}
 
From  Prop.~\ref{GeodEmbedded}, Cor.~\ref{model}.3,  together with the fact that $X$ is defined over $\overline\Q$ whenever $\Gamma$
is arithmetic, we deduce:

\begin{corr} If moreover, $\Gamma\subset G$ is a neat congruence subgroup, and $M$ is embedded in $D$ as a totally geodesic submanifold, then $M$ has a model over a number field.
\end{corr}

\begin{exmples}
1. Since the unit ball $B^n$  in $\C^n$ can be represented  as 
the domain $I_{1,n}$ and since $\ell(I_{1,n})=0$, \emph{all} (positive dimensional) geodesically embedded subvarieties of a compact 
arithmetic quotient of the unit ball have models over a number field.\\
2. A domain of type  $IV_n$ with $n\le 18$  is a parameter space for lattice polarized K3 surfaces, and since $\ell(IV_n)=1$, using
local Torelli,  we deduce that if
we have a family of K3 surfaces over a compact base  $B$ of dimension $\ge 2$  whose  period map is injective and gives a geodesic embdedding, the base manifold 
 $B$ has a model over a number field.
\end{exmples}

\section{Applications to variations of Hodge structure}
\label{sec:appHs}

\subsection{Hodge theory revisited}
\label{ssec:hodgeanew}

As a preliminary to the topic of Shimura varieties,  it is useful to view a Hodge structure  as  a  representation  
space for a certain algebraic torus, as observed by Deligne.  See e.g. \cite[Chap. I]{DMOS}, \cite[Chap. 15]{3authors}.

To explain this briefly, giving  a Hodge structure on a  real vector space $V$ is the same as giving a  morphism 
\[
h:\bS\to \gl V,\quad \bS=\Res_{\C/\R }\G_m,
\]
where I recall that  the Weil restriction $\Res_{\C/\R }\G_m$ is just the group $\C^\times$
considered  as a real group.   In other words, a real Hodge structure  is just a rational (or "algebraic")  representation of the torus group $\bS$.
One sees this by observing that on the complexified vector space
 $V_\C=V\otimes_\R \C$ the action of $\bS$ diagonalizes 
and the Hodge subspace $V^{p,q}\subset V_\C$ by definition  is  the subspace where 
 $h(z)$ acts as multiplication with   $z^p\bar z^q$. 

If  the Hodge structure has pure weight $k$  this shows up  as follows: via the natural inclusion $w:\R^\times \to \bS$, 
the action of $t\in \R^\times$ is multiplication by $t^k$.    This motivates introducing 
\[
w_h= h\comp w: \G_m \to \gl V, \text{the \emph{weight morphism}.}
\]
 If, moreover,   $V$ has a rational structure, say $V=V_\Q\otimes\R$, this
weight morphism  is obviously defined over $\Q$.  When this is the case,
 one defines  the \emph{Mumford-Tate group} of $h$ as  the smallest closed subgroup $M=M(h)$ of $\gl {V_\Q}$ such that
 $h$ factors through the real algebraic group  $M_\R$. 
 
 Hodge structures coming from geometry carry a polarization, where I recall that a \emph{polarization} consists of   a $\Q$-valued bilinear form $b$ on $V_\Q$ 
 satisfying the two Riemann  relations 
\begin{enumerate} 
 	\item $b_\C(x,y) = 0$ if $x$ is in $V^{p,q}$ and $y$
 is in $V^{r,s}$ for $(r,s) \ne (k-p,k-q)$, where $k$ is the weight of the Hodge structure;
	\item $   {\ii^{p-q} }b(x,\overline x) > 0$ if $x$ is a nonzero vector in $V^{p,q}$.  
 \end{enumerate}
A Hodge structure is \emph{polarizable} if such a $b$ exists and then
  $M$ is known to be reductive. See \cite[Prop. I.3.6]{DMOS}, \cite[Prop.~15.2.6]{3authors}.
Using this language, one singles out  a  \emph{CM-Hodge structure} as one  whose  Mumford-Tate group  is abelian and hence, by 
 reductivity, an algebraic torus.
 
 Let me next discuss the notion of a   \emph{variation of Hodge structure}.   It consists of a  local system $\mathbf{V}$ on a smooth quasi-projective  variety $S$  of finite dimensional  $\Q$-vector spaces,
 such that  all fibers admit a polarizable Hodge structure.  More precisely,  $\mathbf{V}$ should come from a representation
of the fundamental group of $S$ in a finite dimensional vector space $V$ equipped with a non-degenerate bilinear form $b$
such that
\begin{enumerate}
\item the locally free sheaf $\VV=V\otimes \OO_S$ carries a descending  filtration $\FF^\bullet$ by holomorphic subbundles;
\item   the natural flat connection $\nabla$ on $\VV$  lowers degrees of this filtration by at most $1$ (Griffiths' transversality);
\item  $b$ and $\FF^\bullet$  induce   a polarized Hodge structure in each stalk.
\end{enumerate} 
  
Given such a variation of Hodge structure, the Hodge structure over $x\in S$ corresponds to $h_x:\bS \to \gl V$ and
its Mumford-Tate group  may vary, However, outside a countable union of proper subvarieties, $M=M(h_x)$ is the same, the \emph{generic Mumford-Tate group}, and
a point with this Mumford-Tate group is called \emph{Hodge generic}.

The group 
\[
G=\aut(V,b)
\]
is  a $\Q$--algebraic group. The 
 representation of  $\pi_1(S,x)$ in $V$ defining the local system $\bV$ preserves the polarization $b$ and the image $\Gamma$ 
 of $\pi_1(S,x)$  in $G(\R)$ is discrete.  It  is called the \emph{monodromy group} of the variation.
\begin{dfn}
The connected component   of the $\Q$--Zariski closure of the monodromy group  in $G$ is called the
\emph{algebraic monodromy group}.
\end{dfn} 
 The group $G(\R)$ acts transitively on a so called \emph{period domain} $D$, which classifies the Hodge structures on $V$
 with a fixed set of Hodge numbers  polarized by $b$. The obvious map $p:S\to \Gamma\backslash D$ is holomorphic; it is called the
 \emph{period map}. The Griffiths' transversality condition is in general a further constraint. It is vacuous for weight one variations
 and also for variations of K3-type.

 \subsection{Application to variations of weight $1$ and $2$} 
 
For a weight two variation with Hodge numbers $h^{2,0}=p, h^{1,1}=q$, the period domain 
has shape  $D= \so{2p,q}{}/\ugr p {}\times \so q{} $, the K3-case corresponding to $p=1,q=19$.
For weight two domains one further introduces the \emph{rank}  $\ell(D)$ of $D$ which generalizes the concept
for  hermitian symmetric spaces from Table~\ref{hsds}:
 \[
 \ell(D)=  \begin{cases} 1 & \text{if  }  p=1,  \\
 q-1 & \text{if  } p=2,\\
 (p-1) t + \epsilon & \text{if  } p\ge 3, t=\lfloor \half(q-1)\rfloor ,\quad\epsilon =
 \begin{cases} 0 & \text{if  }  q \text{ odd}\\
  1& \text{if  } q \text{ even.}
 \end{cases}
\end{cases}
\]
One has the following  rigidity result:
\begin{thm}[\protect{\cite[Theorem 3.1]{rigid}}] \label{thm:rigid} Let $D$ be a period domain for polarized weight 1 or 2 Hodge structures.
An immersive period map  $f: S \to\Gamma \backslash  D$   with $S$ quasi-projective is rigid keeping source and target fixed  as soon as  $\dim S\ge \ell(D)+1$.
\end{thm}

Using Prop.~\ref{GeodEmbedded}, as  a corollary, we get: 
\begin{corr} Let $D$ be a period domain for polarized weight 1 or 2 Hodge structure.  Let   $S$ be quasi-projective and $f: S \to\Gamma \backslash  D$    an immersive period map of rank $\ge \ell(D)+1$.
Suppose moreover, that $S$ is geodesically embedded, then $S$ has a model over $\overline \Q$.
\end{corr}

\subsection{Shimura varieties}
 \label{ssec:shimvars}
 
 One needs a  Hodge theoretic  interpretation of Shimura varieties, i.e., varieties
 the form  $X=\Gamma\backslash D$ for which 
$D=G(\R)/K$ is a Hermitian symmetric domain of non-compact type
and $G$ is a connected $\Q$--algebraic group.   For details of the discussion that follows see e.g.  \cite[Chap. 16,17]{3authors},
\cite{milne}.

A   point $x\in D$ turns out to  correspond  to a unique $h_x:\bS \to G_\R$  and so a given representation of $G$ in $V$
 defines a  real Hodge structure. If  the representation comes from a $\Q$--representation $\rho: G\to \gl {V_\Q}$ one  might not
 get a rational Hodge structure. However, we do get a direct sum of such structures (possibly of different weights) 
 if the weight morphism $   \rho\comp h_x  \comp w:\R^\times \to \gl V$  is defined over $\Q$.
 Such representations exist: take the adjoint representation, with $H=\Lie G$ and $\rho=\rm{ad}:G \to \gl V$: its  weight is zero
 and hence the weight morphism is automatically defined over $\Q$. 
  
 The group $G(\R)$ acts by conjugation on $h_x$. Let $h^{(g)}_x$ denote the conjugate of $h_x$ by $g\in G(\R)$. Then one has the
 basic equality
 \[
 h_{gx}= h_x^{(g)}
 \]
 and hence, since $G(\R)$ acts transitively on $D$, one may view $D$ as an entire conjugacy class of maps $h:\bS\to G_\R$.
 Each  point  in $D$ can be  identified with
 such  a map since $h=h^{(g)}$ precisely if $g$ belongs to the isotropy group of   the corresponding Hodge structure. 
 For clarity, let me write $[h]$ for the point in $D$ corresponding to $h\in \mor(\bS, G(\R))$. 
 Not any $G(\R)$--conjugacy class of a morphism $\bS\to G_\R$ underlies a Hermitian symmetric domain. For this to be true,
such a morphism has to verify certain axioms, as given in \cite{shimuravars}. If this is the case, the corresponding pair $(G,D)$ is called
a \emph{Shimura datum}  and $D$ is called a \emph{Shimura domain}.   By the previous remarks
about the adjoint representation, all Hermitian symmetric domain thus arise with $G$ the group of holomorphic automorphisms of $D$,
which is known to be $\Q$-algebraic and of adjoint type.

 It makes sense to define the  Mumford-Tate group  of a point $[h]\in D$  as the   smallest closed subgroup $M(h)$ of $G$ such that
 $h$ factors through the real algebraic group  $M(h)_\R$.  Then $\rho(M(h))$ is the Mumford-Tate group of the Hodge structure
 $\rho\comp h$. 
 The orbit of $h\in D$ under its Mumford-Tate group  $M(h)$ is   a holomorphic submanifold of $D$ which turns out to be
 a Shimura domain for  $M(h)$. It is called the  \emph{submanifold of Hodge type} 
 passing through $[h]$. Its image in $X$ is called a \emph{special 
 subvariety}. 
 
 As recalled above,  for a  point $[h]\in D$ outside a countable
union of proper closed subvarieties in $D$, the
Mumford-Tate group is precisely $G$.  Call such a point \emph{Hodge-generic}. For such points, $D$ is the submanifold of Hodge type
through $[h]$. At the other end of the spectrum one has the CM-points in $D$, by definition those points $[h]$ for which
 $M(h)$ is abelian (i.e.  an algebraic torus). In this case it is its own submanifold of Hodge type.
Concerning these points, one has:

\begin{lemma}[\protect{\cite[Corr.~17.1.5]{3authors}}] A Shimura subdomain contains a dense set of CM-points.
\end{lemma}

\subsection{Monodromy and rigidity}
The geometry of the variation  is reflected in the \emph{algebraic monodromy}, which as I recall, is the connected component $M^{\rm mon}$ of the
$\Q$--Zariski closure in $\gl V$ of the monodromy group of the variation.  Any reductive group such as $M$ has a canonical almost
direct product decomposition $$M= M^{\rm der} \cdot (\text{center of } M),$$ where $M^{\rm der} $ is the derived subgroup of $M$, its maximal semi-simple
subgroup. 
There are two important results concerning the relation of the two groups: 
\begin{thm} \label{andre}  {\rm 1)\cite[Thm]{andre}} The algebraic monodromy group is a normal subgroup of the generic Mumford-Tate group. In fact, one has $M^{\rm mon}\triangleleft M^{\rm der} $.\\
{\rm 2)\cite[Prop. 2]{andre}}  If there are CM-points in the variation, this is an  equality $M^{\rm mon}= M^{\rm der} $.
\end{thm}

 Let me now consider a more general situation of a homomorphic map $p: S \to \Gamma\backslash D$ to a Shimura variety, i.e.  $D=G(\R)/K$ is a 
 bounded Hermitian symmetric domain. 
This defines a polarizable variation of Hodge structures on $S$ where 
Griffiths' transversality is automatic. Here $\Gamma$ is  the monodromy group of the variation. The group that determines the deformations of  $p$  is the centralizer of the algebraic monodromy group inside the group $G$:
\[
G':= Z_G (M^{\rm mon}).
\]
Indeed, one has:
\begin{prop} Under the assumption that $X=\Gamma\backslash D$ is a  Shimura variety, the  "period map"  $p:S\to \Gamma\backslash D$ is rigid if and only if $G'(\R)$ is compact.
\end{prop} 
\proof The Lie algebra $\eg $ of $G(\R)$ consists of the endomorphisms of $V$ that are skew with respect to $b$. The Cartan involution
induces  a direct sum decomposition $\eg=\ek\oplus\ep$ where $\ek$ is the Lie algebra of the
maximal  compact subgroup  $K(\R)\subset G(\R)$.
The Lie algebra has a natural structure of a weight zero Hodge structure inherited from the one on $\End(V)$. Indeed
\[
\eg_\C= \eg^{-1,1}\oplus \eg^{0,0}\oplus \eg^{1,-1},\quad  \eg^{0,0}= \ek_\C.
\]
The Lie algebra $\eg'\subset \eg$ of $G'(\R)$ consists of those endomorphisms  in $\eg$ that commute with the action of the monodromy group. This subalgebra inherits
a weight zero Hodge structure and
by \cite[Theorem 3.4]{rigid2}, the tangent space to infinitesimal deformations of $p$ is isomorphic to $(\eg'_\C)^{-1,1}$ and in this case, as a real space it is isomorphic to $\ep\cap \eg'$.
Hence $\ep\cap\eg'=0$ if and only if $\eg'=\ek\cap\eg'$ if and only if $G'(\R)$ is compact.
\endproof

Observe next that $G'$ is also a reductive group of Hermitian type: $D_2:=G'(\R)/K\cap G'(\R)$ is a bounded  subdomain of $D$ and  if $\tilde S$ is a universal cover of $S$
with lifting $\tilde p : \tilde S\to D$,
there is an induced holomorphic map $F:S\times D_2 \to D$ extending $\tilde f$. This maps parametrizes the   deformations of $f$ keeping $S$ and $D$ fixed. If  $f$ embeds $\tilde S$   as a subdomain $D_1\subset D$, i.e. if $f$ is a \emph{geodesic} embedding, then one has a product situation
\[
\tilde F: D_1\times D_2\into D.
\]
In other words,  the deformations of the embedding  $\Gamma_1 \backslash D_1 \into \Gamma  \backslash D $ between two Shimura  varieties   
are parametrized by a  Shimura variety of the form $\Gamma_2 \backslash D_2$. By Prop.~\ref{GeodEmbedded}  one  then concludes:
\begin{corr}[\protect{\cite[\S 2]{abdulali}}] \label{CompactandRigid} Let $G$ be a $\Q$-algebraic group of Hermitian type, $G_1\subset G$
a reductive subgroup, and let $D=G(\R)/K$, $D_1=G_1(\R)/K\cap G_1(\R)$ the corresponding domains.  Put  $G_2= Z_G G_1$, $D_2=G_2(\R)/G_2(\R)\cap K$. 
Let  $\Gamma$ be  a neat arithmetic subgroup of $G(\Q)$ such that  $\Gamma_i=\Gamma\cap G_i(\Q)$ $i=1,2$ is  also neat.   The
embedding $\Gamma_1 \backslash D_1 \into \Gamma  \backslash D $ between the corresponding Shimura varieties.
   is rigid with fixed target precisely when $G'(\R)$ is compact.
In particular,   the embedding  $\Gamma_1 \backslash D_1 \times \Gamma_2 \backslash D_2\into \Gamma  \backslash D $  is rigid.
\end{corr}

Let me next consider  the algebraic monodromy group $M^{\rm mon}\subset G$ from an arithmetic perspective. First recall that for any $\Q$-simple algebraic
group $G$  there is a totally real number field $F$ and an absolutely simple  $F$--group $\tilde G$ such that
\[
G=\Res_{F/\Q} \tilde G. 
\]Here $\Res_{F/\Q}$ is the Weil-restriction whereby  an $F$--group is viewed in a functorial way as a $\Q$--group.
For  a  real embedding $\sigma: F \into \R$ let  $\tilde G^\sigma $ be the corresponding conjugate of $\tilde G$. It is called a \emph{factor} of $G$.
Then
\[
G_\R=\prod_{\sigma \in S} \tilde G ^\sigma_\R,\quad S\text{ the set of  embeddings } F\into \R.
\]
Hence, assuming for simplicity that the algebraic monodromy group is simple over $\Q$, one may write:
\[
M^{\rm mon}=  \Res_{F/\Q}\tilde M^{\rm mon} \implies G'= \Res_{F/\Q} Z_{\tilde G} \tilde M^{\rm mon}.
\]
In particular, for every factor  $(\tilde M^{\rm mon})^\sigma$ there is a corresponding  factor  $\tilde G'^\sigma$. This can be
used in the weight one case as follows:

\begin{corr}  \label{NoCompactFactorIsRigid}
Let there be a weight one variation over a quasi-projective variety with $\Q$-simple algebraic monodromy group.
  Assume that $M^{\rm mon}$ has no compact factor.  Then the variation (and the period map) is rigid.
\end{corr}
\proof
In the weight one case, by \cite[Prop. IV.4.3]{satakebook},  $M^{\rm mon}(\R)$ and $G'(\R)$ are in a  sense  "dual": every non-compact factor $(\tilde M^{\rm mon})^\sigma$
corresponds to a compact factor  $\tilde G'^\sigma$. The assumption implies that all factors  of $G'$ must be compact and so the deformation is rigid.
\endproof

This result implies  a quite curious result that states that non-trivial monodromy at the boundary implies rigidity:
\begin{prop}[\protect{\cite[Th. 8.6]{Sa}}]
 A weight one variation over a quasi-projective variety $S$ with a non-trivial unipotent element in
 the monodromy  is rigid.
 
 This holds in particular  if $S$
 is \emph{not} compact and there is at least one non-finite local monodromy operator at infinity.
 
 In these instances,  if  moreover  $S$ is geodesically embedded, it has a model over $\overline\Q$.
 \end{prop}
 \proof 
 First I need a  result about ranks of simple  groups.
  Recall that a reductive $k$--algebraic group $G$ has $k$--rank zero if it has no $k$--split tori.  By  \cite[\S 6.4]{borel} this is 
the case  if and only if $G$ has no non-trivial characters over $k$ and no unipotent elements $g\in G(k)$, $g\not=1$. 
For $k=\R$, the $\R$--rank is zero precisely when  $G_\R$ is compact.

\begin{lemma}  \label{AuxResult} If $G$ is a $\Q$--simple  group  such that $G_\R$ has at least one  compact factor, then the $\Q$--rank of $G$
is zero. In particular, $G$ contains no unipotent elements $g\not=1$.
\end{lemma}
To show this, as before, write  $G=\Res_{F/\Q} \tilde G$ with $\tilde G$ an absolutely simple group defined over a totally
real number field $F$.

A  character $\chi$ for $G$  induces a  character  $\chi_\sigma$ for $\tilde G^\sigma $ and any unipotent $g\in G$ gives a unipotent element $g_\sigma$ in $\tilde G^\sigma$. Suppose $\tilde G^\sigma_\R$ is compact. Then
$\chi_\sigma=1$ and  $g_\sigma=1$ and also $\chi=1$, $g=1$. 
This finishes the proof of the Lemma.
  
The Lemma implies that  the algebraic monodromy group has no compact factors. Hence, by Cor.~\ref{NoCompactFactorIsRigid} the deformation is rigid.
   \endproof

A similar result can be shown for variations of K3-type: 
 
 \begin{prop}[\protect{\cite[Cor. 5.6.3]{SaZ}}]
  Suppose we have a non-isotrivial K3--variation over a quasi-projective variety $S$  with a non-trivial unipotent element in
 the monodromy. Assume that the variation is not isotrivial. Suppose moreover, that
 its rank is \emph{not} 4. Then   the variation (and the period map) is rigid,  and if  $S$ is also  geodesically embedded, then  $S$ has a model over $\overline\Q$.
 \end{prop}
 \proof
 Here  Lemma~\ref{AuxResult}  is used in a different manner.: For a non-isotrivial 
isotypical variation which is non-rigid,  $M^{\rm mon}(\R)$ has one conjugate  isomorphic to $ \slgr{2,\R}$ with representation space $\R^2\otimes\R^2$
and the remaining conjugates  are $\simeq \su{2}$ with representation space $\C^2$.  It follows from the Lemma that the only
possibility to accommodate a non-trivial unipotent element $T$ is when no compact conjugates  are present and  then  the
local system has rank $4$.
\endproof

\subsection{Special subvarieties of Shimura varieties}
\label{ssec:SS}

 Recall (\S \ref{ssec:shimvars}) that a special subvariety of a Shimura variety $X=\Gamma\backslash G/K$, or a subvariety of Hodge type, comes from the orbit of a 
 point  in $D=G/K$ under its own Mumford-Tate group. In this subsection we study them in more detail.
  
  A morphism of Shimura varieties 
 \[
 X_1=\Gamma_1\backslash \, \underbrace{G_1(\R) /K_1 }_{D_1}  \to 
 X_2= \Gamma_2\backslash\, \underbrace{G_2(\R)/K_2}_{D_2}
 \]
  is by definition  induced by an equivariant morphism of Shimura domains.
  Such a morphism is given  by a morphism $\varphi:G_1 \to G_2$ of $\Q$--algebraic groups. It then induces a holomorphic
  maps of Shimura domains $f:D_1\to D_2$  by stipulating that   $f([h_1 ]= [(\varphi\comp h_1)]$ for one
  hence all points $[h_1]\in D_1$.  The Mumford-Tate group of  $[h_1]$ maps under $\varphi$ to the Mumford-Tate group of $f([h_1])$.
 It follows that the subvariety $f(D_1)$ is special in $D_2$: if $[h_1]$ is Hodge generic, then $f([h_1)$ has Mumford-Tate group
 $\varphi(G_1)$ and $f(D_1)$ is the orbit of this group acting on $f([h_1])$.

  Suppose next that $f$ is an embedding. One then may assume that $\varphi$ is also an embedding.
 It is not hard to see that $f$ is a totally  geodesic embedding. See e.g. \cite[Chap 11.5]{3authors}.
 Then the conjugate map $f^{(g')}$, $g'\in G_2(\R)$ is also a totally
 geodesic embedding. It may or may not arise from a morphism of Shimura domains. Indeed, its image may not have
 CM-points at all.
 
 The variation of Hodge structure on $D_2$ restricts to one on  $f^{(g')}D_1$
 and this descends to give one on its image in $X_2$.  The monodromy of this variation is $\Gamma_1^{(g')}$.
 The connected component $M^{\rm mon} $ of its  Zariski closure in $G_2$ acts transitively on $D=f^{(g')}D_1\subset D_2$ and 
 so, if it would be the generic Mumford-Tate group $M$ or, which suffices,
  its derived group, the subdomain $D$ would be   of Hodge type.  
   This is the case if $D$ contains a CM-point. Indeed,  
  by  Theorem~\ref{andre}  one then has $M^{\rm mon}=M^{\rm der}  $. Concluding, I have shown: 
  
 \begin{lemma} Let  $D_1,D_2$ two Shimura domains and let be  $i=f^{(g')}: D_1\into D_2$
  as above.  Then $i$ is a morphism of Shimura domains if and only if  $g'\in G_2(\Q)$ if and only if 
   $i(D_1)$ contains a CM-point of $D_2$.
 \end{lemma}
  
This can be used to give  another proof of Abdulali's criterion  \cite[Thm. 3.4]{abdulali}:

\begin{prop} Let   $i:X_1\into X_2$
be a totally geodesic embedding of Shimura varieties. If the embedding is rigid, $i(X_1)$ is a special subvariety.
\end{prop}
\proof
Since  Shimura varieties are defined over  a number field (cf.~\cite{milne}),  one may  apply  Cor.~\ref{model}.4. So, if the embedding is rigid, the image is defined over a number field.
To show  that the image is a special subvariety, by the previous Lemma,  it suffices to 
 find a CM-point in the image. But,
if $x\in X_1$ is a CM-point, then $i(x)$ is also a CM-point since the Mumford-Tate group
of the Hodge structure corresponding to $x$ is an algebraic torus and hence, so is the one
associated to $i(x)$ since $i$ is defined over $\overline \Q$. 
\endproof

Corollary~\ref{CompactandRigid} then yields examples for weight one Hodge structures:
\begin{exmples}    1. The group  $G_1=\gl 2$ can be embedded in $\smpl g$ as follows.  Set $V_k=(\R^{2k},J_k)$, $J_k=\begin{pmatrix}
0 & \mathbf{1}_k\\
-\mathbf{1}_k & 0
\end{pmatrix}.$  
The direct sum $\oplus_k V_1$ is isomorphic to the symplectic space  $V_k$. Whence  a  faithful representation  
$\rho_k$ of $\sl 2 $:
\[
\begin{pmatrix}
a& b\\
c& d
\end{pmatrix} \mapright{\rho_k} \begin{pmatrix}
a\mathbf{1}_k & b\mathbf{1}_k\\
c \mathbf{1}_k& d\mathbf{1}_k
\end{pmatrix}.
\]
For any $k=1,\dots,g$ the direct sum representation $\rho_k\oplus $ (rank $(g-k)$ trivial representation)
  induces a holomorphic embedding $\eh \into \eh_g$. It  gives the non-compact  embedded Shimura  curves starting from the Shimura datum $(\slgr 2 ,\eh)$.  
There is no locally constant factor if and only if $k=g$ and then  the embedding is rigid.  This follows from Corollary~\ref{NoCompactFactorIsRigid}.
These  non-compact rigid curves are often called  \emph{rigid curves of Satake type}.\\
2. There are also examples where $G_1$ has compact factors. 
Here I use again the Satake "duality" mentioned before, but in its precise form as explained  in \cite[\S 7]{Sa}. It applies to $G_1$ and $G'_1:= Z_{\smpl g} G_1$ and gives:
\[
(G_1)_\R  \simeq \slgr  2 \times \underbrace{\su 2 \times\cdots \times \su 2 }_{m-1} \implies  (G'_1)_\R\simeq \so   2{}   \times \underbrace{\su 1 \times\cdots \times \su 1 }_{m-1}
\]
The latter group is compact and hence the deformation is rigid. There are indeed examples of such embeddings, the \emph{Mumford type curves}. See \cite[\S 6.2]{addington}.
\end{exmples}
%


 %

\section{Applications to Higgs bundles}
\label{sec:Higgs}

 \subsection{Basic notions}
 \label{ssec:basic}
 
 A \emph{Higgs bundle}\footnote{For more details on Higgs bundles see e.g.  \cite[Chapter 13]{3authors}.} over a complex manifold $B$ is a pair $(\VV,\tau)$
 of a holomorphic bundle together with an $\End(\VV)$--valued $1$-form $\tau$ such that $\tau\comp\tau=0$. The form $\tau$ can also be viewed as a  \emph{Higgs field}, a
 homomorphism  $\tau : \VV \to \VV \otimes \Omega^1_B$.  
 A graded   Higgs bundle is a Higgs bundle such that $\VV=\oplus_r  \VV^r$, with $\VV^r$ locally free and such that
 $\tau|\VV^r : \VV^r\to \VV^{r-1}\otimes\Omega^1_B$ .
 
The standard example comes from polarized complex variations of Hodge structures on $B$. Recall \cite[\S 4]
{simpson} that such a system  consist  of
\begin{itemize}
\item   a local system of $\C$--vector spaces 
  $\bV$  equipped with a flat non-degenerate bilinear form. In other words, if $\pi$ is the fundamental group of $B$ based at $o\in S$, $\bV$ comes from a representation $\rho:\pi \to \ogr{V,b} {}$, $V$ the   fiber of $\bV$ at $o$;
  \item a   direct sum decomposition  $\bV\otimes C^\infty_B=\oplus_r \VV_\infty^r$ into locally free $C^\infty_B$--modules
  such that
  \begin{itemize}
                 \item the hermitian form $h(x,y)= (-1)^r b(x,\bar y)$ on $\VV_\infty^r$ is positive definite
                  and the above decomposition is $h$-orthogonal;
                  \item the natural flat connection $\nabla$ on $  \bV\otimes C^\infty_B$ obeys 
  \[
  \begin{matrix}
\quad\VV^r  & \mapright{\nabla}  & \underbrace{\AA^{1,0}_B (\VV_\infty^{r-1)}}_{\downarrow}
		&\oplus&  \underbrace{\AA^1_B(\VV_\infty^r)  }_{\downarrow}&
		 \oplus& \underbrace{\AA^{0,1}_B(\VV_\infty^{r+1})  }_{\downarrow} \\
		&& \tau & + &\quad d &+& \tau^*,
	\end{matrix} 
\]
where $\tau^*$ is the $h$-adjoint of $\tau$.
              \end{itemize}
\end{itemize}
These demands imply that   $\FF^p = \oplus_{r\ge p} \VV_\infty^r$ is a holomorphic subbundle of  $\bV\otimes \OO_B$ and  that Griffiths' transversality holds. This filtration is the \emph{Hodge filtration}.
It also follows that the holomorphic bundle
\[
\VV= \oplus_p  \FF^p/\FF^{p+1},\quad \CC^\infty_B( \FF^p/\FF^{p+1})= \VV^p_\infty
\]
with the  underlying local system  $\bV$ admits the structure of a  graded Higgs bundle with $\tau$  the Higgs field. Flatness (i.e., $ \nabla\comp\nabla=0$) implies the Higgs
condition $\tau\comp\tau=0$. Moreover, the Chern connection, that is, the unique holomorphic connection  on this Higgs bundle  which is metric with respect to
the hermitian metric $h$ turns out to be  $\bar\del +\tau$.
So on any subbundle on which $\tau=\tau^*=0$, the  flat connection $\nabla$ induces  the Chern connection and so the metric $h$
  coincides with the flat metric. Moreover, such a  subbundle comes from a local subsystem of $\bV$ since it is preserved by $ \nabla$.  Also, it is \emph{unitary} since
it admits the flat unitary metric $h$.   This holds in particular for the largest subbundle
for which $\tau=\tau^*=0$:
\[
\mathcal U = \mathbf U \otimes \OO_B :\text{ the \emph{maximal unitary Higgs subbundle}}.
\]
There is an  $h$-orthogonal splitting
\begin{equation}\label{eqn:UniSPlit}
\bV = \mathbf U  \oplus  \bW, \quad  \bW=\mathbf U^\perp.
\end{equation}
 
\subsection{Logarithmic variant}
If $B$ is quasi-projective, one usually considers Higgs bundle with logarithmic growth near the boundary.
To explain this, assume for simplicity that $\dim B=1$ and that $B$ get compactified to a a smooth
 projective curve $\overline B$. Then the boundary $\Sigma=\overline B\setminus B$ consists of finitely many points.
  A  \emph{graded logarithmic Higgs bundle} $\VV=\oplus_p \VV^p$ on $B$, with $\VV^p$ locally free,
  by definition admits  a  Higgs field with components
 \[
\tau: \VV^p\to \VV^{p-1}\otimes \Omega^1_{\overline B}(\log \Sigma).
 \]
  For a variation of Hodge structure on $B$ with unipotent monodromy at the punctures, one lets  $\VV$ be
 the associated graded of the Deligne extended Hodge filtration. Then the Gauss-Manin connection induces a Higgs
 field as above. 
 
 Even more is true.  Choose a coordinate patch $(\Delta,t) $ around a puncture and 
 let $T$ be the (unipotent) local monodromy operator around the puncture. For $v$ a local holomorphic
 section of $\VV$ on the disc, write
 \[
 \nabla v = R \frac{dt}{t}, \quad R\in  \End( \VV|_{\Delta}).
 \]
 Then 
 \[
 N:= R(0) =\log(T) \in \End(V)
 \]
 and the Higgs field at the puncture is given by
\begin{equation}\label{eqn:LogHiggs}
 \tau(0) : \VV_{0} \to \VV_{0}\otimes \frac{dt}{t},\quad  v^p \mapsto (\gr^pN) v^p \otimes \frac{dt}{t}.
 \end{equation}
 
 Suppose  $k$ is the first   index in the grading  for which $\VV^k\not=0$
 and $k+w+1$ the last. Then the number  $w$  is called the width. 
 
 In this general setting, 
 one says that for a Higgs bundle of width $w$,  \emph{the Higgs field is generically maximal}
  if  for all $p\in [k,k+w+1]$ one has
 $\VV^p\not=0$ and if, moreover, $\tau|\VV^p$ generically an is an isomorphism for $p=k,\dots, p+w$.

\subsection{Rigid maximal Higgs subsytems}
\label{ssec:rigidss}

The following   rigidity result  \cite[Lemma 3.1]{carshcurves},  stated without proof, can be formulated in a slightly different way which fits better
within the general framework set up so far:
 \begin{prop} \label{variant2}Let  $B$ be a smooth quasi-projective  variety, 
$\mathbf{V}$ be a local system on $B$ of finite dimensional  $\Q$-vector spaces and let  $\mathbf{W}_\C$ a subsystem
of $\mathbf{V}\otimes\C$.
Suppose $\bW_\C$ is rigid as a subsystem of $\mathbf{V}\otimes\C$. Then $\mathbf{W}_\C$
is  defined over $\overline\Q$ in the sense that $\mathbf{W}_\C= \mathbf{W}\otimes\C$, where
$\bW$ is a local system   of   $F$-vector spaces for  some number field $F$.
\end{prop}
\proof
Let $\pi$ be the fundamental group of $B(\C)$ based at $o\in B(\C)$ and let $V$ be the fiber at $o$ of  $\bV_\C$,  considered
as a $\pi$-representation space.   
The group $\pi$ acts on the  Grassmannian  $G(r,V)$ of $r$-dimensional subspaces $W\subset V$, where $r=\rank \mathbf{W}_\C$.
A fixed point $[W]$ of this  action
corresponds to a complex subsystem of $\bV\otimes\C$. More precisely, the corresponding
$\pi$-invariant subspace is     the fiber $U_{[W]}$ at $[W]$ of the tautological bundle $U\to G(r,V)$.   The spread of the point $[W]$  is a subvariety $Y\subset G(r,V)$ contained in  the locus of fixed points of the $\pi$-action,
because  this action is defined over $\Q$.
The tautological subbundle over $Y$  gives a deformation of $\mathbf{W}_\C\subset \mathbf{V}_\C$, and so rigidity implies that  
  $Y(\overline \Q)=[W]$,  
an isolated point. Hence the local system $\mathbf{W}_\C$, which  corresponds to $U_{[W]}$, is defined over $\overline \Q$.
  \endproof

\medskip

As a direct application, one has:
 \begin{prop}[\protect{\cite[Lemma 3.3]{carshcurves} }] \label{zuovw} Let $B=\overline B\setminus \Sigma$ as above.
 Suppose that $\bV$ underlies a polarizable  $\Q$-variation of Hodge structure and let $(\VV,\sigma)$ the
 corresponding graded logarithmic Higgsbundle over $ B$ with unipotent  monodromy around points of $\Sigma$.
 Let $\bV_\C= \mathbf U  \oplus  \bW$ the splitting  \eqref{eqn:UniSPlit}.
Suppose that the logarithmic Higgs subbundle   $\mathcal W$   corresponding  to $\bW$ is a generically maximal  Higgs subbundle.  Then the  splitting  is  defined over $\overline\Q$. 
 \end{prop}
\proof 
To show how this result  is implied by Proposition~\ref{variant2},  it is enough to show that $ \bW$ is rigidly embedded in  $\bV$.
Again, with $V$ a typical fiber of $\bV$, small deformations of   $ \bW$ are parametrized by the tangent space  to  the fixed locus  
under the $\pi$-action on  the Grassmannian $G(r,V)$  at a $\pi$-invariant point $[W]$.
A tangent vector is therefore represented by
a homomorphism of local systems
\[
q:  \bW \to   \bV/ \bW= \mathbf U
\]
which  is compatible with the  structure as a complex system of Hodge bundles: a small deformation of $\bW$ within $\bV$ inherits
this structure from the one on $\bV$ and the map $q$ is the embedding of the deformed $\bV$ followed by restriction to $\mathbf U$.
But   the Higgs field   for  the left hand is generically an isomorphism
while on the right hand it is zero.  This is impossible unless $q=0$. 
\endproof

\begin{rmq}
A variant of this (loc. cit.) is when $\bW$ is a direct sum of complex systems of Hodge bundles of
\emph{different} widths, all with generically maximal Higgs field. Then almost the same argument
shows that also this splitting is defined over $\overline\Q$.
There is one subtlety here: one has to compare projections  between complex systems of different widths and then one needs  semi-simplicity 
for variations of Hodge structures. This property  is a highly non-trivial consequence of another rigidity property
due to Schmid \cite{vhs}.
See \cite[\S 5]{mono} for details.
 \end{rmq}
 
 One can say more, if  there are punctures with infinite local monodromy:
 \begin{prop} \label{UniPot} The situation is as  in Prop.~\ref{zuovw}. In particular, all local monodromy operators at the 
 boundary are unipotent. Assume  that at least one local monodromy operator has infinite order.
 Then the splitting $
\bV = \mathbf U  \oplus  \bW
$ from \eqref{eqn:UniSPlit} is defined over $\Q$ and $\mathbf U$ extends as a local system to $\overline B$. 
The monodromy of this last system is finite.
 \end{prop}
 \proof
 The property that a Higgs field is an isomorphism on $B$
 extends to  $\overline B$. If all  graded fields $\tau_p$ are isomorphisms, at a puncture \eqref{eqn:LogHiggs} implies  that the $\gr^p N$ are
 isomorphisms and hence that $N$ is an isomorphism. This holds for $\mathcal W$, while  the fact that the Higgs field for $\mathcal
 U$ remains
 zero at  a puncture implies that the local monodromy for $\mathbf U$ is the identity and thus that this local system extends
 to $\overline B$.

 Suppose that we have a splitting as above, valid over a Galois extension $F/\Q$. The property that $N$ is an isomorphism or zero is
 preserved by the action of the Galois group $\GG$. It follows that for $\sigma\in\GG$ 
the  natural inclusion followed by projection
 \[
 \mathbf U^\sigma  \to   \bV \to \bV/ \mathbf U= \bW
\]
sends the fiber  $ \mathbf U_s^\sigma$ at a puncture  $s\in\Sigma$ 
to zero. In other words $ \mathbf U_s= \mathbf U_s^\sigma$    for all $\sigma\in\GG$ and so
this fiber is defined over $\Q$.  Since $\mathbf U$ extends to $\overline B$,
and since the entire
monodromy action is  defined over $\Q$, the local system $\bU$ which is built from the monodromy representation
on some fiber $\bU_s$ is then defined over $\Q$.  The polarization is defined over $\Q$ as well and hence $\mathbf W=\bU^\perp$ is defined over $\Q$. The finiteness of the monodromy follows since the system is defined over $\Q$ and the polarization $h$  on it
is a positive definite  $\Q$-valued form   preserved  by the monodromy.
\endproof 

\begin{exmple}[An interesting Shimura curve in the Torelli locus]  The above result definitely fails when $B=\overline B$: the global monodromy
of $\bU$ may be infinite.  The simplest example  from  \cite[Example 5.1]{Torloc}  is 
a Shimura curve and can be described as follows.
Consider the family of projective curves with affine equation
\[
y^5= x(x-1)(x-t).
\]
This gives a family $C_t$ over $\bP^1$ of genus $4$ curves. The fibers are smooth for $t\not= 0,1,\infty$.
Note however that local monodromy operators are quasi-unipotent in this case, but this does not really
matter since this could be taken care of by  a finite branched cover of $\bP^1$. For simplicity this will not be done since
the above  analysis still works after some minor modification.

Let $\zeta_5$ be primitive root of unity. Then the cyclic group $\Z/5\Z$ generated by
$(x,y)\mapsto (\zeta_5y,x)$ preserves $C_t$ and the  Hodge structure $H^1(C_t)$ of  weight one admit an action of $\Z/5\Z$. 
Let $F=\Q(\zeta_5)$. The Galois group $\GG$ of $F/\Q$ is generated by the element $\sigma$ which sends $\zeta_5$ to $\zeta_5^2$.
It permutes the eigenspaces of $\Z/5\Z$ acting on $V_t=H^1(C_t,\C)$ as in the following table.
\begin{table}[ht]
\caption{Eigenspaces for $\Z/5\Z$ on $V_t$ }
\begin{center}
\begin{tabular}{|c|c|c|}
\hline
Eigenvalue& $h^{1,0} $& $h^{0,1}$ \\
\hline
$\zeta_5$& $0$ & $2$ \\
\hline
$\sigma(\zeta_5)=\zeta_5^2$& $2$& $0$ \\
\hline
$\sigma^2(\zeta_5)= \zeta_5^4$&$1$&$1$\\
\hline
$\sigma^3(\zeta_5)=\zeta_5^3$&$1$&$1$\\
\hline
\end{tabular}
\end{center}
\end{table}
Next, consider the splitting of the corresponding Higgs bundle. 
The Higgs bundle  splits also in eigenspace subbundles;
the Higgs field is zero for the subbundles corresponding to   the  first two rows and an isomorphism for those corresponding to the last two. So the Higgs bundle is maximal and 
one has
\[
\bU = V_{\zeta_5}\oplus V_{\zeta_5^2},\quad \bV=  V_{\zeta_5^4}\oplus V_{\zeta_5^3}.
\]
If  $M^{\rm mon}$ is the algebraic monodromy group of this family, one has a corresponding splitting
\[
M^{\rm mon}(\R)= \su{2} \times \su{1,1}.
\]
The actual monodromy group $\Gamma$ in this case is dense in both factors and hence cannot be finite for the 
local system $\bU$. Proposition~\ref{UniPot} then implies
that the \emph{local monodromy} around the punctures cannot be of infinite order. Indeed, one can
show that the local monodromy operators are
all of order $5$ in this case.  Hence the period map extends over the punctures and the period map embeds the base curve as a compact curve in the period domain. This is consistent with the fact that
fibers over the punctures are  of so called compact type: their generalized Jacobians are products of principally polarized
Abelian varieties whose dimensions sum up to $4$. 

Note also that this is an  elementary  example  giving a negative answer to the following question of Fujita \cite{katata}: 

"for a family $f:X\to B$  of complex
algebraic manifolds over a curve $B$, is  the sheaf $f_*\omega_{X/B}$ is semi-ample?"

 In the above example the latter  sheaf is just
the graded part $\mathcal H^{1,0}= \mathcal U ^{1,0}\oplus \mathcal V^{1,0}$ of the Higgs bundle and while the second bundle is
ample, the first is flat and would be semi-ample if and only if its \emph{global} monodromy would be finite, but this  group is dense in $M^{\rm mon}(\R)$. Compare also the much more elaborate
examples in \cite{fuj}.

\end{exmple}

\end{document}